\documentclass{gtpart} 

\agtart

\usepackage{appendix}
\usepackage{mathrsfs}
\usepackage[usenames]{color}
\usepackage[all, dvips]{xy}
\usepackage{verbatim,rotating,tikz}
\usepackage{setspace}

\title{Cosimplicial models for the limit of the Goodwillie tower}

\author{Rosona Eldred}
\givenname{Rosona}
\surname{Eldred}
\address{Fachbereich Mathematik der Universit\"{a}t Hamburg\\
D- 20146 Hamburg\\
Germany
}
\email{rosona.eldred@math.uni-hamburg.de}
\urladdr{http://www.math.uni-hamburg.de/home/eldred/}


%Homotopy Functors 
%Localization and completion 
%Homotopy equivalences
\arxivreference{math.AT/1108.0114} 
\arxivpassword{ 2fgwj}

\volumenumber{} 
\issuenumber{}
\publicationyear{} 
\papernumber{} 
\startpage{} 
\endpage{}
\doi{} 
\MR{} 
\Zbl{} 
\received{} 
\revised{} 
\accepted{} 
\published{} 
\publishedonline{} 
\proposed{} 
\seconded{} 
\corresponding{} 
\editor{} 
\version{}

\theoremstyle{plain}

\newtheorem{thm}{Theorem}[section]
\newtheorem*{thm*}{Theorem}
\newtheorem{lem}[thm]{Lemma}
\newtheorem{defn}[thm]{Definition}
\newtheorem*{defn*}{Definition}
\newtheorem{prop}[thm]{Proposition}
\newtheorem{cor}[thm]{Corollary}

\theoremstyle{definition}
\newtheorem{rem}[thm]{Remark}

\newcommand{\la}{\leftarrow}
\newcommand{\ra}{\rightarrow}
\newcommand{\lra}{\longrightarrow}
\newcommand{\Mor}{\text{Mor}}
\newcommand{\Top}{\mathrm{Top}}
\newcommand{\sSet}{\mathrm{sSet}}
\DeclareMathOperator{\Hom}{Hom}
\newcommand{\Wedge}{\vee}
\newcommand{\join}{\ast}

\newcommand{\ox}{\otimes}
\newcommand{\x}{\times}

\newcommand{\pt}{\bullet}
\newcommand{\DDelta}{\Delta^{\pt}}
\newcommand{\PP}{\mathrm{P}}
\renewcommand{\P}{\mathrm{P}}
\newcommand{\vDelta}{\Delta}%{\bigtriangleup}
\newcommand{\diag}{\mathrm{diag}}
\newcommand{\Sing}{\mathrm{Sing}}
\newcommand{\T}{\mathrm{T}}

\DeclareMathOperator{\sk}{sk}
\DeclareMathOperator{\cosk}{cosk}
\DeclareMathOperator{\holim}{holim}

\DeclareMathOperator{\hocolim}{hocolim}

\newcommand{\I}{\mathbb{I}}
\newcommand{\X}{\mathscr{X}}
\newcommand{\Y}{\mathscr{Y}}
\renewcommand{\C}{\mathscr{C}}

\begin{document}
%%%%%%%%%%%%%%%%%%%%%%%%%%%%%%%%%%%%%% 

%%%%%%%%%%%%%%%%%%%%%%%%%%%%%%%%%%%%%% 
\begin{abstract}
We call attention to the intermediate constructions $\T_n F$ in Goodwillie's Calculus of homotopy functors, giving a new model which naturally gives rise to a family of towers filtering the Taylor tower of a functor. We also establish a surprising equivalence between the homotopy inverse limits of these towers and the homotopy inverse limits of certain cosimplicial resolutions. This equivalence gives a greatly simplified construction for the homotopy inverse limit of the Taylor tower of a functor $F$ under general assumptions. 
\end{abstract}

%%%%%%%%%%%%%%%%%%%%%%%%%%%%%%%%%%%%%% 
\maketitle

\section{Introduction}

Let $\Delta^n$ be the $n$-simplex and $\sk_0 \Delta^n$ be its $0$-skeleton, that is, $n+1$ points. We use $\ast$ to denote the topological join. For a space $X$, $\sk_0 \Delta^0 \join X \sim CX$ and $\sk_0 \Delta^n \join X \sim \bigvee_n \Sigma X$. Thus, we have the cosimplicial space 
\[
(\sk_0 \Delta^\ast \join X) 
\sim 
\xymatrix{
CX   \ar@<5 pt>[r]\ar@<-5 pt>[r]& 
\Sigma X\ar[l]  \ar@<8 pt>[r] \ar[r] \ar@<-8 pt>[r]& 
\Sigma X  \Wedge \Sigma X  \ar@<-4pt>[l]\ar@<4pt>[l] \ar@<8 pt>[r]\ar@<-5 pt>[r]^-{\vdots}& 
\cdots
}
\] 
%\smallskip
%
%
Hopkins \cite{hopkins, hopkinsthesis}, and later, Goerss \cite{Goerss-Barratt}, analyzed the spectral sequence associated to this cosimplicial space and showed that when $X$ is connected, it converges to $\Z_\infty X$, the Bousfield $\Z$-nilpotent completion of a space (for details of the construction and applications, see Bousfield-Kan\cite{BK}). 

One result of this paper is a new conceptual proof of this theorem. We assume for the moment that $X$ is a connected space. $\P_\infty \I (X)$ is the inverse limit of the Goodwillie Taylor tower of the identity functor, applied to $X$. We show the following weak equivalence, where $\I$ is the identity functor of spaces: 
\[
\holim \sk_0 \Delta^\ast \join X \sim \PP_\infty \I (X).
\]  
By work of Arone-Kankaanrinta  \cite{AroneK-Id}, $\PP_\infty \I (X) \sim \Z_\infty (X)$. This gives the result of Goerss and of Hopkins. 

Indeed, we show that $\forall k \geq 0$ and $X$ connected, that $\holim \sk_k \Delta^\ast \join X \sim \PP_\infty \I (X)$.

More generally, if $F$ is $\rho$-analytic (see Section \ref{sec:bkgd}\; for a definition) and we denote by conn$(X)$ the connectivity of $X$, % is at least $m \geq -2$, % connected for $m \geq -1$, 
%$-1 \leq k= conn(X) \leq \rho$,  
then we have weak equivalences,
\[
\PP_\infty F (X) \sim \holim_\Delta F(\sk_k \Delta^\ast \join X),
\]
for all $k \geq max(\rho-conn(X)-1, 0)$.% \rho-m$ if $\rho \geq m$, and otherwise it is true for all $k\geq 0$. 

This arises as a natural corollary of our main results (specifically Corollary \ref{cor:rho}, following from Theorems \ref{thm:main} and \ref{thm:Tnk} described below).
In order to precisely state our main results, we require some additional definitions and notation. \\

Let $F$ be a functor from spaces to spaces (or to spectra) which preserves weak equivalences. We will also assume that $F$ commutes with filtered colimits.   
Goodwillie, in \cite{GC1,GC2,GC3}, constructs for such $F$ a tower of functors which approximate $F(X)$ under mild conditions on $F$ and its input, $X$. Each finite stage of the tower is denoted $\PP_n F$, in analogy with the $n$th partial sums of a Taylor series of a function, which are degree $n$ polynomials. A full definition of a degree $n$ polynomial functor is given in Section \ref{sec:bkgd}. The homotopy inverse limit of the Taylor tower is denoted $\PP_\infty F$.   %

In general, to build each of the $\PP_nF$'s requires taking the homotopy colimit over a directed system. Each finite stage of this system is the iteration of a homotopy inverse limit construction, called $\T_n F$.  %We will define here the construction, reserving further explanation until Section \ref{sec:bkgd}. 

Let $\mathscr{P}([n])$ be the power set  on $[n] \in \Delta$, a poset. Posets may be viewed naturally as categories with maps given by the partial order, i.e. inclusion of subsets. We denote by $\mathscr{P}_0([n])$ the full subcategory omitting the empty set.

We then define  
%\[
$\T_n F(X) := \holim_{U \in \mathscr{P}_0([n])} F(U \join X)$
%\]
 and observe that this construction gives us a natural transformation $t_n: F(X) \ra \T_n F(X)$, since $F(X) = F(\emptyset \join X)$ and the inclusion of the empty set into each $U$ (viewed as an element of $\mathscr{P}{([n])}$) is compatible with the other maps in the diagram. The $t_n$'s give rise to the maps in the directed system used to construct $\PP_n F(X)$: 

\[
\PP_n F(X) := \hocolim_k (\T_n F(X) \overset{t_n}{ \lra} \T_n^2 F(X) \overset{t_n}{ \lra} \cdots  \overset{t_n}{ \lra}\T_n^k F(X) \overset{t_n}{ \lra} \cdots)
\] 

%%%%
\smallskip 
%%%%

%Following a suggestion of Tom Goodwillie, w
We first establish a new model for each iterated approximation, $\T^k_n F$. Using this model, we obtain maps $\tau^k: \T_{n+1}^k F(X) \ra \T_n^k F(X)$, and therefore a tower of partial approximations for each $k$: 
 \[
\cdots \T^k_n F(X) \overset{\tau^k}{\ra} \T^k_{n-1} F(X)  \overset{\tau^k}{\ra} \cdots \T^k_1 F(X)
\]

%\noindent We will now give our model for the $\T_n^k F$'s, which will be used in the proof of our main theorem:
\noindent We will need a few definitions, and then may give our model for the $\T_n^k F$'s, which will be used in the proof of our main theorem.

Let $\Delta$ be the category of finite ordered sets and monotone maps, with objects denoted $[j]=\{0,\ldots, j\}$ with the usual order.  Let $\Delta_{\leq n}$ be the full subcategory of $\Delta$ with objects $[j]$ such that $j \leq n$. We use diag to denote the diagonal of a $k$-cosimplicial space. 
%

%%%%%%%%%%%%
The $n$-coskeleton of $X^\pt$, denoted $\cosk_n X^\pt$,  is constructed by precomposing with the inclusion $\Delta_{\leq n} \hookrightarrow \Delta$ and then taking the right Kan extension along the inclusion of the subcategory. We let $\cosk_{\vec{n}}(X^{\vec{\bullet}})$ denote the $k$-cosimplicial analog, $n$-coskeleton taken in every dimension.
%\fi

%%%%%%%%%%%%%%%%%%%%%%%%%%%%%%%%% 
%%%%%%%%%%%%%%%%%%%%%%%%%%%%%%%%%
\begin{thm} \label{thm:Tnk} 
%There is a weak equivalence 
For all $k, n \geq 0$, we have the following weak equivalence
\[
\T^k_n F(X) \sim \holim_{\Delta{\leq nk}} \diag ( \cosk_{\vec{n}} F((\sk_0 \DDelta)^{\join k} \join X))
\]
In particular, as $n \ra \infty$, we have as an immediate consequence the following equivalence:  %
\[
\holim_n \T_n^k F(X) \sim \holim_\Delta \diag F((\sk_0 \DDelta)^{\join k} \join X)).
\]
\end{thm}
%%%%%%%%%%%%%%%%%%%%%%%%%%%%%%%%%
%%%%%%%%%%%%%%%%%%%%%%%%%%%%%%%%%

\noindent  We also show that 

%%%%%%%%%%%%%%%%%
\begin{thm}\label{thm:main}
For all $k \geq 0$, the functors $\sk_k \DDelta$ and $(\sk_0 \DDelta)^{\join (k+1)}$ are both homotopy left cofinal as functors from $\Delta$ to $(k-1)$-connected spaces of CW type; in particular, for all spaces $X$ and %nonempty spaces $X$ and 
homotopy endofunctors $F$, we have weak equivalences
\[
\holim_\Delta F(\sk_k \DDelta \join X)  \sim \holim_\Delta  F( (\sk_0 \DDelta)^{\join(k+1)}\join X)
\]
Furthermore, with Theorem \ref{thm:Tnk}, we have weak equivalences for all $k \geq 0$
\[
\holim_{\Delta} F(\sk_k \DDelta \join X) %
\sim \holim_n (\cdots \ra \T^{k+1}_n F(X) \ra \T^{k+1}_{n-1} F(X) \ra \cdots \T^{k+1}_1 F(X))
\]
\end{thm}
%%%%%%%%%%%%%%%%%

The weak equivalences in Theorems \ref{thm:Tnk} and \ref{thm:main} are natural in $k$. We discuss the definition and main property of left cofinal functors in Section \ref{sec:bkgd}.

%We will discuss the maps involved in these towers in section \ref{sec:bkgd}. %, as we make use of lemmas from that section to make the maps clearer.

Notice that the tower in Theorem \ref{thm:main} is over the maps $\tau^{k+1}: \T_n^{k+1}F \ra \T_{n-1}^{k+1}F$, along the same stage of iteration of \textit{different} $\T_n F$'s. This is markedly different than the directed system used to construct the $\PP_n F$'s, which is over the maps $t_n: \T_n^j F \ra \T_n^{j+1}F$, i.e. along iterations of the \textit{same} $\T_n$ construction. 
  We depict both collections of maps in Figure \ref{fig:part-tower}.  $\PP_n F$ is then the homotopy colimit along the $n$th column of Figure \ref{fig:part-tower}, whereas the partial approximation towers in Theorem \ref{thm:main} are the rows.

\begin{figure}[h]
\[
%\scalebox{.8}{$
\begin{array}{ccccccccc|cc}          
\P_\infty F := \holim & (\cdots \lra  &\P_n F& \lra  & \P_{n-1} F         & \lra&          &  \lra& \P_1 F)&       & \\
&                & \underset{\rotatebox{270}{(}}{\rotatebox{270}{$= \hocolim$}} & &\underset{\rotatebox{270}{(}}{\rotatebox{270}{$= \hocolim$}} & & &  &\underset{\rotatebox{270}{(}}{\rotatebox{270}{$= \hocolim$}}&\\                
&                & \uparrow  &      & \uparrow          &      &          &      &  \uparrow  &  \\
&                & \vdots &      & \vdots         &      &          &      &  \vdots &       &\\                
&                & \uparrow  &      & \uparrow          &      &          &      &  \uparrow  &  &\text{row:}\\                
&\cdots \lra & \T_n^2 F & \lra & \T_{n-1}^2 F & \lra & \cdots & \lra & \T_1^2 F& &    1\\
&                & \uparrow  &      & \uparrow          &      &          &      &  \uparrow &  \\
&\cdots \lra & \underset{\rotatebox{270}{)}}{\T_n F} & \lra & \underset{\rotatebox{270}{)}}{\T_{n-1}} F & \lra & \cdots & \lra & \underset{\rotatebox{270}{)}}{\T_1 F} & & 0 \\ \hline      
&\text{column:} & n && (n-1)&&&&1& & \\          
\end{array}
%$}
\]
\caption{Partial approximations}
\label{fig:part-tower}
\end{figure}

%\noindent  
Theorem \ref{thm:main} provides an equivalence between the homotopy limit over the $k$th row of the diagram in Figure \ref{fig:part-tower} and the homotopy limit of $F$ applied to $(\sk_k \Delta^\ast \join X)$.  We will discuss the maps involved in these towers in Section \ref{sec:bkgd}.

%%%%
\bigskip
%%%%

\noindent We now present several consequences of Theorem \ref{thm:main}.

%%%%%%%%%%%%%%%%%%%%%%%%%%%
%
% [list corollaries]
%
%%%%%%%%%%%%%%%%%%%%%%%%%%%

%%%%%%%%%
\begin{cor}\label{cor:bejm}\footnote{This can be seen as the unstable extension of a stable pro-result in Bauer-Eldred-Johnson-McCarthy \cite[Prop 4.2]{un}}
For a given $j$, there are weak equivalences for all $k\geq 0$
\[
\PP_j F(X) \sim \holim_\Delta \PP_j F(\sk_k \DDelta \join X)
\]
In particular, this also implies that $P_\infty F(X) \sim \holim_\Delta P_\infty F(\sk_k \DDelta \join X)$.
\end{cor}
%%%%%%%%%

\noindent A functor $F$ is said to be $\rho$-analytic if its failure to be polynomial of any degree is bounded in a way that depends on $\rho$. One important consequence of $\rho$-analyticity is that for $X$ at least $\rho$-connected, $F(X) \sim \P_\infty F(X)$. This is the consequence of analyticity which we need for our purposes; for a thorough definition of analyticity and discussion of its implications, see  \cite{GC2}.

We will use this fact and the preceding corollary to establish the following:

%%%%%%%%%
\begin{cor}\label{cor:rho}
Let $F$ be a $\rho$-analytic functor. %Let $X$ be a space that is nonempty, that is, at least $-1$-connected. 
Then we have weak equivalences $\forall  k > \rho$, 
\[
\scalebox{.9}{$
\P_\infty F (X) \sim \holim_n ( \cdots \T^{k+1}_n F(X) \overset{\tau^{k+1}}{\ra} \T^{k+1}_{n-1} F(X)  \overset{\tau^{k+1}}{\ra} \cdots \overset{\tau^{k+1}}{\ra}\T^{k+1}_1 F(X)) \sim \holim_\Delta F(\sk_k \Delta^\ast \join X)$}
\]
\noindent If we raise the connectivity of $X$ (denoted conn(X)), we may improve this to all $k \geq \text{max} (\rho-\text{conn}(X)-1, 0)$.

\end{cor}
%%%%%%%%%

\begin{proof}
For two spaces, $X$ and $Y$, the connectivity of their join, $X \join Y$, is equal to $\text{conn}(X)+\text{conn}(Y)+2$. For each $m$, $\sk_k \Delta^m$ is $(k-1)$-connected, which yields  conn($\sk_k \Delta^m \join X) =(k+1)+\text{conn}(X)$.  We conclude that the connectivity of the cosimplicial space $\sk_k \DDelta \join X$ is the same.  In particular,  $X$ arbitrary (of connectivity $\geq -2$), $\sk_k \DDelta \join X$ has connectivity at least $(k-1)$ and for $F$ of analyticity $\rho$ such that $\rho-1 \leq k$, $F(\sk_k \DDelta \join X) \sim \P_\infty F(\sk_k \DDelta \join X)$.

We apply holim to this equality, and also recall that by Cor \ref{cor:bejm}, we have $\P_\infty F(X) \sim \holim_\Delta \P_\infty F(\sk_k \DDelta \join X)$. The resulting equivalence is  
%\[
$\P_\infty F(X) \sim \holim_\Delta F(\sk_k \DDelta \join X).$
%\]
We follow with Theorem \ref{thm:main} and obtain the other stated equivalence. 
\end{proof}
\smallskip

We note that this corollary yields $\PP_\infty F(X)$ via a construction which no longer makes use of any \textit{infinite} colimits, and that each $\T_n^k$ commutes with holim and hofiber. 

%%%
\smallskip
%%%

%\noindent 
If $F$ is $\rho$-analytic, and $X$ is in its radius of convergence (i.e. at least $\rho$-connected), then Cor \ref{cor:rho}\; implies
\[
F(X) \sim \PP_\infty F(X) \sim \holim_\Delta F(\sk_0 \DDelta \join X). 
\]

\noindent This result may be rephrased as saying that for $F$ an analytic functor, and a space $X$ in its ''radius of convergence'', $F (X)$ is well approximated by $F(- \;\join X)$ applied to finite nonempty sets.

 By a similar proof to that showing that the identity is 1-analytic\footnote{For the identity functor, this may be found in \cite[Example 4.3, Theorem 2.3]{GC2}}, it follows that if a functor $F$ commutes with realizations and preserves filtered colimits, then it is 1-analytic.   Corollary \ref{cor:rho} then gives us the following: 
 
%%%%%%%%%
\begin{cor}If $F$ commutes with realizations and preserves filtered colimits, the equivalence
\[
\PP_{\infty} F(X) \sim \holim_n ( \cdots \T^{k+1}_n F(X) \overset{\tau^{k+1}}{\ra} \T^{k+1}_{n-1} F(X)  \overset{\tau^{k+1}}{\ra} \cdots \T^{k+1}_1 F(X)).
\]
holds for all $k>1$.
\end{cor}

\begin{prop}\label{cor:Znil1}
For the identity $\I$ from spaces to spaces, $X$ a connected space, and $\Z_\infty X$ the Bousfield $\Z$-nilpotent completion of $X$, we have that for all $k \geq 0$, the following weak equivalence
\[
\holim_\Delta (\sk_k \DDelta \join X ) \sim \PP_\infty (\I) (X) \sim \Z_\infty X,
\] 
and when $X$ is already nilpotent, $\holim_\Delta (\sk_k \DDelta \join X )  \sim X$. 
\end{prop}

\begin{proof}
Since $\I$ is 1-analytic, Corollary \ref{cor:rho} allows us to conclude that for any space $X$, $\holim_\Delta (\sk_k \DDelta \join X) \sim \PP_\infty (\I) (X) $ for all $k\geq 2$. Restricting to $X$ which are 0-connected, i.e. raising the minimum connectivity from $-2$ to 0, changes this equivalence to hold for all $k\geq 0$.  Then, by Arone-Kankaanrinta \cite[\S 3]{AroneK-Id}, we have that for $X$ connected, $ \PP_\infty (\I) (X) \sim \Z_\infty X$, and for $X$ nilpotent, $\PP_\infty (\I) (X) \sim X$. 
\end{proof}

%%%%
\smallskip
%%%%

If, instead of combining our result with that of Arone-Kankaanrinta, we combine with the result of Goerss and Hopkins, one can view our main result as justification for why the spectral sequences associated to the Taylor Tower of the identity of spaces and that associated to the $\Z$-nilpotent completion of a space abut to the same thing \footnote{The spectral sequence associated to the Taylor Tower of $\I(X)$ takes as input the collection of $D_n \I (X)$, which were computed by Brenda Johnson \cite{BrendaThesis}.}.

%%%%
\smallskip
%%%%

%%% Added Remark, ''Futurework/applications''
\begin{rem}
In highlighting these partial approximation towers, we also raise the question of what  being weakly equivalent to $\P_\infty F(X)$ off of the radius of convergence is telling us. 

This should in general provide conditions that give the ``true'' radius of convergence of a functor. That is, given $F$ a $\rho$-analytic functor, it will also converge for spaces $X$ such that $F(X) \sim \holim_\Delta F(\sk_k \DDelta \join X)$ for all $k \geq \text{min}(\rho-m-1, 0)$, where $X$ is at least $m$-connected. The example we would like to use is slightly outside of the scope of this paper as written, but follows from the appropriately general versions of the same theorems stated in the dissertation of the author \cite{me}.

Given a map $f:A \ra B$ of commutative ``rings'' (simplicial, or spectra), we can then tensor (or smash) $B$ over $A$ with a set, and by prolongation, with a cosimplicial simplicial set by tensoring levelwise. That is, we can consider cosimplicial rings $\sk_q \DDelta \ox_A B$. We then either work with augmented $A$-algebras or with the unbased calculus of Bauer-Johnson-McCarthy\cite{un2} with $f$ as our relevant map that we are factoring. In both settings, the identity functor is 1-analytic, and we have that $\P_\infty \I (B) \sim \holim \sk_1 \DDelta \ox_A B$.  If the map $f$ is flat and \'{e}tale, then $\holim \sk_1 \DDelta \ox_A B \simeq B$, with no connectedness assumptions on the map $f$, and $B \simeq \P_\infty \I(B)$. That is, we have the extension of the radius of convergence then to such $B$. 
%We believe there should be a way to interpret \cite[Thm 3.2.2.]{DGM} as an instance of this.

We also suspect that the $k$th partial approximation tower, $\holim_\Delta F(\sk_k \DDelta \join X)$, is in general a sort of $k$-analytic approximation of a $\rho$-analytic functor, where $k$ is less than $\rho$. Evidence for this can be seen in the case of the identity functor, which is 1-analytic. For $X$ 0-connected, the 0th partial approximation tower is equivalent to $\Z_\infty (X)$, and thanks to \cite[p31]{AroneK-Id}, we know that $\Z_\infty (X)$ is 0-analytic.
\end{rem}

%%%%
\bigskip
%%%%

\noindent \textit{Note on translating between the language of Goerss \cite{Goerss-Barratt} and our current terminology:} %\\
Goerss defines a cosimplicial construction $C(X,X)$ for a nonempty space $X$ which sends $[n]$ to $\bigvee_{n} \Sigma X$  (taking the empty wedge here as $CX$). As mentioned earlier, this is equivalent to $(\sk_0 \Delta^\ast \join X)$. 
Theorem 1.1 in \cite{Goerss-Barratt} implies that when $X$ is connected, $\holim_\Delta C(X,X) \sim \Z_\infty X$, that is,  $\holim_n \T_n \I (X):= \holim_\Delta (\sk_0 \Delta^\ast \join X) \sim \Z_\infty X$. 
\bigskip

\noindent \textit{Note on translating between the language of Hopkins and our current terminology: }\\
\noindent This is from \cite{hopkins}, Section 3, p221-222.  He lets $C_n$ be what we call $\mathscr{P}_0 ([n])$. He defines, for a given space $X$, a functor $F^n$, as the homotopy inverse limit of a (punctured) cube. For $A \in \mathscr{P}_0 ([n])=:C_n$,  the $A$-indexed position of this $(n+1)$-cube is the homotopy colimit of X mapping to $|A|$ different copies of a point, which we will explain shortly. He denotes this by $F^n A$. Regarding $A$ as a finite ordered set, we can view $F^nA$ as the homotopy pushout of the following: 

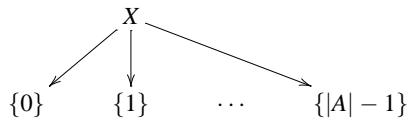
\[
\scalebox{.85}{$
\xymatrix{
                 &X\ar[dl]\ar[d]\ar[drr]\\
\{0\} & \{1\} & \cdots & \{|A|-1\} \\
} 
$}
\]

\noindent We replace these maps by cofibrations (since we are taking a homotopy colimit), giving us that we are pushing out over the following diagram: 

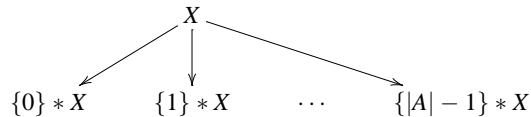
\[
\scalebox{.85}{$
\xymatrix{
                 &X\ar[dl]\ar[d]\ar[drr]\\
\{0\}\join X & \{1\}\join X & \cdots & \{|A|-1\}\join X \\
} 
$}
\]

\noindent That is, the $A$-indexed position of this $(n+1)$ cube is  $A\join X$. 

Then $F^n := \holim_{A \in  \mathscr{P}_0 ([n])} F^n A \sim\holim_{U \mathscr{P}_0 ([n])} U \join X$. That is, we have shown that his $F^n$'s exactly the $\T_n \I(X)$'s.  He constructs a tower of these $F^n$'s: 
 
\[
(\cdots \ra \holim F^n \ra \holim F^{n-1} \ra \cdots \holim F^1)
\]
which is therefore our $\T_n \I$ tower, 
 \[
 (\cdots \ra \T_n\I(X) \ra \T_{n-1}\I(X) \ra \cdots \T_1 \I(X)).
 \]
 
 \noindent Theorem 3.2.2 of \cite{hopkinsthesis} is that the homotopy inverse limit of a construction that is equivalent to the  tower of $F^n$'s gives $\Z_\infty X$ when $X$ is connected, i.e.  that $\holim_n \T_n \I(X) \sim \Z_\infty X$.

%%%%%%%%%%%%%%%%%%%%%%%%%%%%%%%%%%%%%% 
%%%%%%%%%%%%%%%%%%%%%%%%%%%%%%%%%%%%%%
%
% ORGANIZATION (SUB)section 
\subsection{Organization} 
%
%
%%%%%%%%%%%%%%%%%%%%%%%%%%%%%%%%%%%%%%
%%%%%%%%%%%%%%%%%%%%%%%%%%%%%%%%%%%%%%

The remainder of this paper is organized as follows. 
Section \ref{sec:bkgd} gives background mainly on the %
Calculus. Theorem \ref{thm:Tnk}\;, the more geometric and cosimplicial interpretation of $\T_n F$, is proven in Section \ref{sec:Tnk}. The proof of Theorem \ref{thm:main} is then given in Section \ref{sec:main}.

%%%%%%%%%%%%%%%%%%%%%%%%%%%%%%%%%%%%%% 
%%%%%%%%%%%%%%%%%%%%%%%%%%%%%%%%%%%%%% 
%
%
% ACKNOWLEDGEMENTS
\subsection{Acknowledgements}
%
%
%%%%%%%%%%%%%%%%%%%%%%%%%%%%%%%%%%%%%% 
%%%%%%%%%%%%%%%%%%%%%%%%%%%%%%%%%%%%%% 

The author would like to thank Tom Goodwillie for giving generously of his time to discussions that led to the conjectural form of the main results of this paper. The observations and work that led to that point would not have been possible without the patience and guidance of Randy McCarthy and the current forms of the proofs benefitted greatly by discussions with him. A discussion with Bill Dwyer was very helpful for streamlining one of the proofs, as well as gaining more insight into cofinality arguments.  Thanks also go to Victor Turchin, who pointed out an error in an earlier statement of Theorem \ref{thm:Tnk}. 
%%%
The author also gratefully acknowledges the help of the people who gave feedback on earlier drafts of this work, including Bert Guillou, Kristine Bauer and Brian Munson, and the helpful comments of the reviewer. 

This work was partially supported by travel funds from the Midwest Topology Network (NSF grant DMS-0844249).

%%%%%%%%%%%%%%%%%%%%%%%%%%%%%%%%%%%%%%
%%%%%%%%%%%%%%%%%%%%%%%%%%%%%%%%%%%%%%
%
%  BACKGROUND section
\section{Background}% on Goodwillie Calculus}
\label{sec:bkgd} 
%
%
%%%%%%%%%%%%%%%%%%%%%%%%%%%%%%%%%%%%%%
%%%%%%%%%%%%%%%%%%%%%%%%%%%%%%%%%%%%%%
%\subsection{Goodwillie Calculus}
We will first discuss cofinality, and then review the required constructions of Goodwillie Calculus. 
For this paper, we restrict our attention to functors $F$ from spaces to spaces (not necessarily based) which commute with filtered colimits.  $F$ is a \textit{homotopy functor} if it preserves weak equivalences.

%%%%%%%%%%%%%%%%%%%%%%%%%%%%%%%%%%%%%%
%%%%%%%%%%%%%%%%%%%%%%%%%%%%%%%%%%%%%%
\subsection{Homotopy left cofinality}
%%%%%%%%%%%%%%%%%%%%%%%%%%%%%%%%%%%%%%
%%%%%%%%%%%%%%%%%%%%%%%%%%%%%%%%%%%%%%
This is the definition of homotopy left cofinal which we will use. There are also corresponding notions of non-homotopy cofinality (involving strict limits) and a dual notion of right cofinality (involving colimits) which we will not discuss. 

%%%%%%%%%
\begin{defn}[see Hirschhorn \cite{hirschhorn}, Definition 19.6.1 p418] \label{def:cof}
\label{def:cof2}
Let $\mathscr{D}$ be a small category and for all $A,B$ objects of $\mathscr{D}$, we denote by $\Mor_{\mathscr{D}}(A,B)$ the set of morphisms in $\mathscr{D}$ between them. Let $G: \Delta \ra \mathscr{D}$.  The functor $G$ is \textbf{homotopy left cofinal} if for every object $\alpha$ of $\mathscr{D}$, the simplicial set $n \mapsto \Mor_{\mathscr{D}} (G(n), \alpha)$ is contractible.
  \end{defn}
%%%%%%%%%

The following consequence of being homotopy left cofinal is what we use to establish our equivalences: 

%%%%%%%%%
\begin{lem}[see Hirschhorn\cite{hirschhorn}, Theorem 19.6.7 \& 16.6.23]\label{thm:hirsch}
Let $M$ be a simplicial model category, let $\mathscr{C, D}$ be small categories, and let $G: \mathscr{C} \ra \mathscr{D} $ be a functor.  If $G$ is \textbf{homotopy left cofinal}, then for every object-wise fibrant $\mathscr{D}$-diagram $F$ in $M$, we have that the following natural map of homotopy limits is a weak equivalence:
\[
\holim_{\mathscr{D}} F \ra \holim_{\mathscr{C}} F \circ G
\]
\end{lem}
%%%%%%%%%

We will be working only with a simplicial model category $M$ where the objects are all fibrant, namely, the category of spaces with cofibrations the cellular inclusions of CW complexes, so object-wise fibrancy comes for free. 

Next, we will review the Calculus.

%%%%%%%%%%%%%%%%%%%%%%%%%%%%%%%%%%%%%%
%%%%%%%%%%%%%%%%%%%%%%%%%%%%%%%%%%%%%%
\subsection{Excisive functors}
%%%%%%%%%%%%%%%%%%%%%%%%%%%%%%%%%%%%%%
%%%%%%%%%%%%%%%%%%%%%%%%%%%%%%%%%%%%%%
In \cite{GC1}, Goodwillie establishes the following condition for a functor, which is in analogy with a function being polynomial of degree 1: 

\begin{defn}
A functor is \emph{1-excisive} if $F$ takes homotopy pushout (called \emph{cocartesian}) squares to homotopy pullback (called \emph{cartesian}) squares. 
\end{defn}

\noindent This may not be the most familiar statement of excision, compared to its usual statement as one of the axioms of a (generalized) homology theory (as in Eilenberg-Steenrod\cite{ES}) . There is a nice discussion in the dissertation of Mauer-Oates\cite[p.22]{MO-Thesis} of how to get from excision as usually stated in the Eilenberg-Steenrod axioms to this definition.

%%%%%%%%%%%%%%%%%%%%%%%%%%%%%%%%%%%%%%
%%%%%%%%%%%%%%%%%%%%%%%%%%%%%%%%%%%%%%
\subsection{Excisive approximation} 
%%%%%%%%%%%%%%%%%%%%%%%%%%%%%%%%%%%%%%
%%%%%%%%%%%%%%%%%%%%%%%%%%%%%%%%%%%%%%

%

The following is a pushout of finite sets. It is also a diagrammatic representation of the category $\mathscr{P} ([1])$. 
\[
\xymatrix{
\emptyset \ar[r]\ar[d] & \{0\}\ar[d]\\
\{1\} \ar[r] & \{0, 1\}\\
}
\]
We make the following definition:
\[
\T_1 F(X) := \holim_{U \in \mathscr{P}_0 ([1])} F(U \join X)
\]
As a result,  we have a natural transformation $t_1: F(X) \ra \T_1 F(X)$, induced by the natural map: 
\[
F(X)=F(\emptyset \join X) \ra \underset{U \in \mathscr{P}_0 ([1])}{\holim} ( U \mapsto F(U \join X)).
\]
That is, the map from the initial object of the square, $F(X)$ to the homotopy pullback of the rest, $\T_1 F(X)$. We can take $\T_1$ of $\T_1 F$, and also have the same natural transformation from initial to homotopy pullback, now $t_1: \T_1 F(X) \ra \T_1 (\T_1 F(X)) =: \T_1^2 F(X)$. See Figure \ref{fig:T_1^2}. 

We define the 1-excisive approximation to $F$, $\P_1 F$, as the following homotopy colimit over iterations:
\[
\P_1 F(X) := \hocolim ( \T_1 F(X) \overset{t_1}{\ra} \T_1^2 F(X) \overset{t_1}{\ra} \cdots)
\]

%%%%%%%%%%%%%%%%%%%%%%%%%%%%%%%%%%%%%%
%\subsection{Iterate:}
%%%%%%%%%%%%%%%%%%%%%%%%%%%%%%%%%%%%%%

%%%%%%%%%%%%%%%%%%%%%%%%%%%%%%
\begin{figure}[h!]
\[
\scalebox{.8}{$
\begin{array}{rcl}
%%%%%%%%%%%
% FirstEquality
%%%%%%%%%%%
\T_1^2 F(X) &:= &   \holim %_{\mathscr{P}_0([1])} 
\left( 
	\begin{array}{c}
	\xymatrix{
	 	& \T_1 F(\{0\}\join X) \ar[d]\\
	 \T_1 F(\{1\} \join X) \ar[r] & \T_1 F(\{0,1\} \join X)\\
	}
	\end{array}
\right)
\\
\\
\\
%%%%%%%%%%%
% Third Equality
%%%%%%%%%%%
&\sim&  %\underset{\mathscr{P}_0([1])\x \mathscr{P}_0([1])} {\holim}
\holim 
\left(
	\begin{array}{ccc}
	& &%%BLOCK%%
		\scalebox{.75}{$\left( 
		\begin{array}{c}
		\xymatrix{
								&F(\{0\}\join \{0\} \join X)\ar[d]\\
		F(\{0\}\join \{1\} \join X) \ar[r] & F(\{0\}\join\{0,1\} \join X)\\
		}
		\end{array}
		\right)$}
	%%END BLOCK%%
		\\
		&&\downarrow\\
	%%BLOCK%%	
	\scalebox{.75}{$\left( 
		\begin{array}{c}
		\xymatrix{
								&F(\{1\}\join \{0\} \join X)\ar[d]\\
		F(\{1\}\join \{1\} \join X) \ar[r] & F(\{1\}\join\{0,1\} \join X)\\
		}
		\end{array}
		\right)$}
	%%END BLOCK%%	
		&\rightarrow&
	%%BLOCK%%	
	\scalebox{.75}{$\left( 
		\begin{array}{c}
		\xymatrix{
								&F(\{0,1\}\join \{0\} \join X)\ar[d]\\
		F(\{0,1\}\join \{1\} \join X) \ar[r] & F(\{0,1\}\join\{0,1\} \join X)\\
		}
		\end{array}
		\right)$}
	%%END BLOCK%%	
		\\		
	\end{array}
	\right)
\end{array}
$}
\]
\caption{$\T_1^2F(X)$}
\label{fig:T_1^2}
\end{figure}
%%%%%%%%%%%%%%%%%%%%%%%%%%%%%%

%%%%%%%%%%%%%%%%%%%%%%%%%%%%%%%%%%%%%%
%%%%%%%%%%%%%%%%%%%%%%%%%%%%%%%%%%%%%%
\subsection{Higher Degree Functors}
%%%%%%%%%%%%%%%%%%%%%%%%%%%%%%%%%%%%%%
%%%%%%%%%%%%%%%%%%%%%%%%%%%%%%%%%%%%%%
As in the 1-excisive case, we begin with a diagrammatical representation of the powerset category, now $\mathscr{P}([n])$, which is an $(n+1)$-cube indexed by subsets of $[n]$. 

%%%%%%%%%%%
\begin{defn}
We say that a $\mathscr{P}([n])$-indexed diagram (i.e. an $(n+1)$-cube) $\X$ is \emph{strongly co-cartesian} if every square (i.e. 2-dimensional) sub-face is cocartesian.  
\end{defn}

\begin{defn}
We say that a functor $F$ is \emph{$n$-excisive} if it takes strongly co-cartesian $(n+1)$-cubes to cartesian $(n+1)$ cubes.
\end{defn}

Analogous to the 1-excisive case, we make the following definition:
\[
\T_n F(X) := \holim_{U \in \mathscr{P}_0 ([n])} F(U \join X)
\]
This allows us to express $t_n: F(X) \ra \T_n F(X)$ as the natural map:
\[
F(X)=F(\emptyset \join X) \ra \underset{U \in \mathscr{P}_0 ([n])}{\holim} ( U \mapsto F(U \join X)).
\]
As before, we define the degree $n$ polynomial approximation to $F$, $\P_n F$, as the following homotopy colimit:
\[
\P_n F(X) := \hocolim ( \T_n F(X) \overset{t_n}{\ra} \T_n^2 F(X) \overset{t_n}{\ra} \cdots)
\]
%%%%%%%%%%%%%%%%%%%%%%%%%%%%%%%%%%%%%%
%%%%%%%%%%%%%%%%%%%%%%%%%%%%%%%%%%%%%%
\subsection{Taylor Tower}
%%%%%%%%%%%%%%%%%%%%%%%%%%%%%%%%%%%%%%
%%%%%%%%%%%%%%%%%%%%%%%%%%%%%%%%%%%%%% 
The collection of polynomial approximations to a functor $F, \{ \P_n F\}_{n \geq 0}$, come with natural maps $\P_nF(X) \ra \P_{n-1} F(X)$ for all $n \geq 1$. Using Goodwillie's \cite{GC3} definition
 \[
 (\T_n^iF)(X) := \underset{(U_1, \ldots, U_i) \in \mathscr{P}_0 ([n+1])^i}{\holim} F(X \join (U_1 \join \cdots \join U_i ))
\]
we then have for all $i,n \geq 1$ a natural map $\T_n^i F \overset{q_{n,i}}{\longrightarrow} \T_{n-1}^i F$ induced by the inclusion of categories, $\mathscr{P}_0 ([n])^i \hookrightarrow \mathscr{P}_0 ([n+1])^i$. Taking the colimit along $i$ gives us the induced map $\P_n F \overset{q_{n}}{\longrightarrow}\P_{n-1} F$. With these maps we form a tower, the Goodwillie Taylor Tower of $F(X)$:
\[
\cdots \ra \P_n F(X)  \overset{q_{n}}{\longrightarrow} \P_{n-1} F(X) \overset{q_{n-1}}{\longrightarrow} \cdots  \longrightarrow \P_1 F(X) \overset{q_{1}}{\longrightarrow} \P_0 F(X)
\]

\noindent We denote by $\P_\infty F(X)$ the homotopy inverse limit of this tower. 
%\smallskip

%It is worth noting now that t
This defines $\P_\infty F(X)$ as the %
homotopy inverse limit of a collection of constructions which are themselves % 
homotopy colimits (of finite homotopy inverse limits). That is, it is not expected that this construction will commonly commute with either colimits or limits; there are several special cases set out in \cite{GC3}. 

%%%%%%%%%%%%%%%%%%%%%%%%%%%%%%%%
\subsection{Analyticity} 
%%%%%%%%%%%%%%%%%%%%%%%%%%%%%%%%
Let $\rho$ be an integer greater than or equal to zero. We say that a functor $F$ is \textit{$\rho$-analytic} if its failure to be $n$-excisive is controlled by $\rho$ as $n$ increases. A precise definition of $\rho$-analyticity may be found in \cite{GC2}. We need only the following consequence of being $\rho$-analytic.

For a $\rho$-analytic functor, $F$, and any $\rho$-connected space, $X$, the following natural map is an equivalence $F(X)\overset{\sim}{\ra}\PP_\infty F(X)$.  \\

\noindent \textbf{Note:} Higher values of $\rho$ mean that $X$ is `closer' to 0 (i.e. $\ast$), since increasing connectivity means that $X$ has more vanishing homotopy groups. A \textit{lower} value of $\rho$ means a larger ``radius of convergence'' of the functor $F$. 

%%%%%%%%%%%%%%%%%%%%%%%%%%%%%%%%%%%%%%
%%%%%%%%%%%%%%%%%%%%%%%%%%%%%%%%%%%%%%
%%%%%%%%%%%%%%%%%%%%%%%%%%%%%%%%%%%%%%
%%%%%%%%%%%%%%%%%%%%%%%%%%%%%%%%%%%%%%
%
% Proof of Ancillary therem
\section{Proof of Theorem \ref{thm:Tnk}}\label{sec:Tnk}
%
%%%%%%%%%%%%%%%%%%%%%%%%%%%%%%%%%%%%%%
%%%%%%%%%%%%%%%%%%%%%%%%%%%%%%%%%%%%%%

\noindent In this section, we prove the following theorem:\\

%%%%%%%%%%%%%%%%%%%%%%%%%%%%%%%%%%%%
%%%%%%%%%%%%%%%%%%%%%%%%%%%%%%%%%%%%
\noindent \textbf{Theorem \ref{thm:Tnk}} \textit{For all $k, n \geq 0$, we have the following weak equivalence}
\[
\T^k_n F(X) \sim \holim_{\Delta{\leq nk}} \diag ( \cosk_{\vec{n}} F((\sk_0 \DDelta)^{\join k} \join X))
\]
\textit{In particular, as $n \ra \infty$, we have as an immediate consequence the following equivalence: } %
\[
\holim_n \T_n^k F(X) \sim \holim_\Delta \diag F((\sk_0 \DDelta)^{\join k} \join X)).
\]
%%%%%%%%%%%%%%%%%%%%%%%%%%%%%%%%%
%%%%%%%%%%%%%%%%%%%%%%%%%%%%%%%%%
\smallskip 

\noindent The heart of this proof is establishing a finite-dimensional analog for the generalized $n$-cosimplicial Eilenberg-Zilber Theorem (see Shipley \cite[Prop 8.1]{shipcos}).% which says that for $X \; n$-cosimplicial, we have an equivalence $\holim_\Delta \diag (X)\sim \holim_{\Delta \x \cdots \x \Delta} X$. 

First, we make use of a lemma to switch models for $\T_n F$. %, all $n$. 

%%%%%%%%%%%%%%%%%%%%%%%%%%%%%%%%%
\begin{lem}[Sinha \cite{Dev} Theorem 6.7, or Hopkins \cite{hopkinsthesis}, \S 3.1 Prop 3.1.4] \label{prop:G_n}
Let $c_n:\mathscr{P}_\emptyset ([n])\ra \Delta_{\leq n}$ be the functor which sends a nonempty subset subset $S$ to $[\# S-1]$ and which sends an inclusion $S\subseteq S'$ to the composite $[\#S-1]\cong S \subset S'\cong [\#S'-1]$.  $c_n$ is homotopy left cofinal. 
\end{lem}
%%%%%%%%%%%%%%%%%%%%%%%%%%%%%%%%%

The immediate consequence of this lemma (plus Lemma \ref{thm:hirsch}, which outlines the relevant consequence of cofinality) is that we can move between the two models for each $\T_n F$, i.e. the following are weakly equivalent:
\[
\holim_{\Delta_{\leq n}} F(\sk_0 \DDelta \join X) \overset{\sim}{\longleftarrow} \holim_{U \subset \mathscr{P}_0([n])} F(U \join X)=: \T_n F(X)
\]

For the next step of the proof, we will need to introduce some new notation. Consider $X$ an $n$-cosimplicial space and, for various $i$ and $j$, $\cosk_j^i X$ an $n$-cosimplicial space with the $j$-coskeleton taken in the $i$th direction. This will be explained for bicosimplicial objects, and the constructions generalize  to the $n$-cosimplicial case.  
\medskip

Let $X^{\bullet, \bullet}$ be bicosimplicial, with directions 1 and 2. For each $p$ and each $q$, we have cosimplicial spaces $X^{p, \bullet}$ and $X^{\bullet, q}$. For each $p$, we can consider (separately, not necessarily for all $p$ at once) the cosimplicial space $\cosk_j X^{p, \bullet}$ (and similarly for each $q$). Then, the bicosimplicial space $X$ with $j$-coskeleton taken in the 1st direction, denoted $\cosk_j^1 X$, is the functor $q \mapsto (p \mapsto  \cosk_j X^{p, \bullet})^q$. 

We also note that these then give rise to partial holims %or totalizations 
taken in the various directions. Let $X$ be an $n$-cosimplicial space, then the homotopy limit taken only in the $i$th direction, $\holim_{\Delta}^i X$ and the $j$th holim in the $i$th direction, $\holim_{\Delta_{\leq j}}^i X$, both produce $(n-1)$-cosimplicial spaces. 

%%%%%%%%%%%%%%%%%%%%%%%%%%%%%%%%%
\begin{prop}\label{stablem}
For any $n$-cosimplicial space, $Y$, we have that 
\[
\holim^1_{\Delta_{\leq j_1}}\holim^2_{\Delta_{\leq j_2}} \cdots\holim^n_{\Delta_{\leq j_n}} Y \sim \holim_{\Delta \leq J} \diag \cosk_{j_1,j_2,\ldots, j_n} Y,
\]
%i.e. $\Tot^1_{j_1}\Tot^2_{j_2} \cdots \Tot^n_{j_n} Y \sim \Tot_{J} \diag \cosk_{j_1,j_2,\ldots, j_n} Y$ for $Y$ fibrant. 
Here, $\cosk_{j_1,j_2,\ldots, j_n}$ denotes the $j_i$-th coskeleton taken in the $i$th direction for all $i \in \{1,\ldots, n\}$.
\end{prop}
%%%%%%%%%%%%%%%%%%%%%%%%%%%%%%%%%

In the special case that $X$ be an $n$-cosimplicial space with the property that 
$\holim_{\Delta}^i X$, is equivalent to  $\holim_{\Delta_{\leq j_i}}^i X$, for each $i \in \{1,\ldots, n\}$ and some $j_i  \geq 1$, this is Lemma 2.9 of Bauer-Eldred-Johnson-McCarthy \cite{un}, albeit stated in slightly different language and using an alternate method of proof. 
\smallskip

We first establish the bicosimplicial case, and then conclude $n$-cosimplicial by induction.  Bicosimplicial spaces can be assumed to be objectwise fibrant, so $\holim_\Delta X$ can be written $\Hom (\Delta, X)$, which will simplify the proof.

%%%%%%%%%%%%%%%%%%%%%%%%%%%%%%%%%%%%
%%%%%%%%%%%%%%%%%%%%%%%%%%%%%%%%%%%%
\begin{lem}\label{bico}
Let $X$ be a bicosimplicial simplicial set. Then 
\[
\holim_p^1 \holim_q^2 X = \holim (\cosk_p^1 \cosk_q^2 X) \sim \holim_{p+q} \diag (\cosk_p^1 \cosk_q^2 X)
\]
\end{lem}
%%%%%%%%%%%%%%%%%%%%%%%%%%%%%%%%%%%%
%%%%%%%%%%%%%%%%%%%%%%%%%%%%%%%%%%%%

Let sSets  be the category of simplicial sets, then $sSets^{\Delta \x \Delta}$ is the category of bicosimplicial simplicial sets. 

%%%%%%%%%%%%%%%%%%%%%%%%%%%%%%%%%%%%
\begin{proof}[Proof of Lemma \ref{bico}]
%%%%%%%%%%%%%%%%%%%%%%%%%%%%%%%%%%%%
We will freely exchange $\holim (X)$ with its model $\Hom_{sSets^{\Delta \x \Delta}}(\Delta^{s_1}_t \x \Delta^{s_2}_t , X^{s_1,s_2}_t)$ (enriched Hom). Using Proposition 8.1 of Shipley\cite{shipcos}, this is also equivalent to $\Hom_{sSets^{\Delta}} (\Delta^s_t , X^{s,s}_t)$. Likewise, 
%\[
%%%%%%%%%
%\medskip
\begin{equation}\label{blah}
\holim (\cosk_p^1 \cosk_q^2 X) \cong \Hom_{sSet^{\Delta}}(\Delta^s_t, (\cosk_p^1 \cosk_q^2 X)^{s,s}_t).
\end{equation}
\medskip
%%%%%%%%%
%\]
\noindent By the Yoneda Lemma applied to $X \in sSets^{\Delta \x \Delta}$, 
\[
\begin{aligned}
a,b,c \mapsto  X^{a,b}_c &\cong a,b,c \mapsto \Hom_{Set^{\Delta \x \Delta}} (\Delta_a^\pt \x \Delta_b^\ast, X^{\pt,\ast}_c) \\
a,c \mapsto X^{a,a}_c 
&\cong  a,c \mapsto \Hom_{Set^{\Delta }} (\Delta_a^\pt, \diag (X^{\pt, \ast})_c)\\
\end{aligned}
\]
\smallskip

\noindent Given this, and combined with the isomorphisms $\holim_n (X) \cong \Hom (\sk_n \Delta, X) \cong \Hom (\Delta, \cosk_n X)$,  we also have 
\[\begin{aligned}
a,b,c \mapsto  (\cosk^1_p \cosk^2_q X)^{a,b}_c &\cong a,b,c \mapsto \Hom_{Set^{\Delta \x \Delta}} ((\sk_p \Delta^\pt)_a \x (\sk_q \Delta^\ast)_b, X^{\pt,\ast}_c). 
\end{aligned}
\]
%
%%%%%%%%%
%\noindent This allows us to then provide the 
We continue the proof, making use of this in the first of the following equivalences: % to (\ref{blah}),
\begin{spacing}{1.5}
\[
\scalebox{.85}{$
\begin{array}{rcll}
(\ref{blah}) &\cong& \Hom_{sSet^{\Delta}}(\Delta^s_t, \Hom_{Set^{\Delta \x \Delta}} (\sk_p \DDelta_s \x \sk_q \Delta^\ast_s, X^{\ast, \pt}_t))\\
%\\
&\cong& \Hom_{sSet^{\Delta \x \Delta}}(\Delta^s_t \ox_\Delta (\sk_p \DDelta_s \x \sk_q \Delta^\ast_s), X^{\ast, \pt}_t) & \Hom_{Set^{\Delta \x \Delta}} (G, -\;) \text{ right adjoint to }-\ox_\Delta G \\
%\\
&\cong& \Hom_{sSet^{\Delta \x \Delta}}(\Delta^s_t \ox_\Delta \sk_{p+q}(\sk_p \DDelta_s \x \sk_q \Delta^\ast_s), X^{\ast, \pt}_t) &(\sk_p \DDelta_s \x \sk_q \Delta^\ast_s) \cong \sk_{p+q} (\sk_p \DDelta_s \x \sk_q \Delta^\ast_s)\\
%\\
&\cong& \Hom_{sSet^{\Delta \x \Delta}}(\sk_{p+q}\Delta^s_t \ox_\Delta (\sk_p \DDelta_s \x \sk_q \Delta^\ast_s), X^{\ast, \pt}_t)& \sk_n \text{ is a tensor,}  \sk_nY = \Delta\otimes_{\Delta_{\leq n}} Y\\
%\\
& \cong& \Hom_{sSet^{\Delta}}(\sk_{p+q}\Delta^s_t,  \Hom_{Set^{\Delta \x \Delta}} (\sk_p \DDelta_s \x \sk_q \Delta^\ast_s, X^{\ast, \pt}_t) & \Hom, \ox \text{ adjunction again}\\
%\\
& \cong&  \Hom_{sSet^{\Delta}}(\sk_{p+q}\Delta^s_t, (\cosk_p^1 \cosk_q^2 X)^{s,s}_t) & \text{Definition/Yoneda }\\
%\\
& \cong & \holim_{p+q} \diag (\cosk_p^1 \cosk_q^2 X) \\
\end{array}$}
\]
\end{spacing}

We made use of the fact that $n$-skeleton is a tensor, by means of the following yoga:%. This is important, as 
\[
\Delta^*\otimes \sk_nY = \Delta^*\otimes_{\Delta} \Delta \otimes_{\Delta_{\leq n}} Y = \Delta^*\otimes_{\Delta_{\leq n}} Y = 
 \Delta^*\otimes_{\Delta_{\leq n}} \otimes_{\Delta}\Delta \otimes _{\Delta} Y = sk_n\Delta^*\otimes_{\Delta} Y
 \]
%%%%%%%%%%%%%%%%%%%%%%%%%%%%%%%%%%%%
\end{proof}
%%%%%%%%%%%%%%%%%%%%%%%%%%%%%%%%%%%%
%%%

%%%%%%%%%%%%%%%%%%%%%%%%%%%%%%%%%%%%
\begin{proof}[Proof of  Prop \ref{stablem}]
%%%%%%%%%%%%%%%%%%%%%%%%%%%%%%%%%%%%
%Recall that $X$ is $n$-cosimplicial and $\Tot^i_{j_i} X \sim \Tot^i X$ for all $i \in \{1, \dots, n\}$ . 
Our inductive hypothesis is that 
\[
\holim_{j_{n-1}}^{n-1}  \cdots  \holim_{j_2}^2  \holim_{j_1}^1X \sim \holim_{\sum \limits_{i=1}^{n-1}j_i} \diag X.
\]
\noindent The diagonal is taken in the directions we are totalizing in, which we will abusively denote as $\diag X$; it should be clear from context what is meant.

Now consider $\holim_{j_n}^n \holim_{j_{n-1}}^{n-1}  \cdots  \holim_{j_2}^2  \holim_{j_1}^1X$. 
\[
\begin{array}{rcll}
\holim_{j_n}^n \holim_{j_{n-1}}^{n-1}  \cdots  \holim_{j_2}^2  \holim_{j_1}^1X &\sim&  \holim_{j_n} \holim_{\sum \limits_{i=1}^{n-1}j_i} \diag X & \text{Inductive Hypothesis}\\
&\sim&  \holim_{\sum \limits_{i=1}^{n}j_i} \diag (\diag X )& \text{Lemma }\ref{bico}\\
&\sim&  \holim_{\sum \limits_{i=1}^{n}j_i} \diag X \\
\end{array}
\]
The last step is that taking the diagonal of the entire object is equivalent to taking the diagonal in the first $k$ dimensions and \textit{then} the diagonal of that collapsed out part and the rest.
%
%%%%%%%%%%%%%%%%%%%%%%%%%%%%%%%%%%%%
\end{proof}
%%%%%%%%%%%%%%%%%%%%%%%%%%%%%%%%%%%%

%%%%%%%%%%%%%%%%%%%%%%%%%%%%%%%%%%%%
\begin{proof}[Proof of Theorem \ref{thm:Tnk}]
%%%%%%%%%%%%%%%%%%%%%%%%%%%%%%%%%%%%
Now that we have Proposition \ref{stablem}, this is nearly immediate. 
\begin{spacing}{1.5}
\[
\begin{array}{rcll}
\T_n^{k} F(X) &:= & \T_n (\T_n ( \cdots (\T_n F(X))\\
			   & \sim& \holim_{\Delta_{\leq n}}\cdots \holim_{\Delta_{\leq n}}  F( (\sk_0 \DDelta)^{\join k} \join X) & \text{Lemma \ref{prop:G_n}}\\
			   & \sim& \holim_{\Delta_{\leq n}\x \cdots \Delta_{\leq n}} F( (\sk_0 \DDelta)^{\join k} \join X) & \text{Fubini}\\
			  & \sim& \holim_{\Delta_{\leq nk}} \diag \cosk_{\vec{n}} F( (\sk_0 \DDelta)^{\join k} \join X) & \text{Prop \ref{stablem}}%\text{Cosimplicial Eilenberg-Zilber}\\
			  
\end{array}
\]
\end{spacing}

%%%%%%%%%%%%%%%%%%%%%%%%%%%%%%%%%%%%
\end{proof}
%%%%%%%%%%%%%%%%%%%%%%%%%%%%%%%%%%%%

%%%%%%%%%%%%%%%%%%%%%%%%%%%%%%%%%%%%%% 
%%%%%%%%%%%%%%%%%%%%%%%%%%%%%%%%%%%%%% 
%
%
\section{Proof of Theorem \ref{thm:main}}\label{sec:main}
%
%
%%%%%%%%%%%%%%%%%%%%%%%%%%%%%%%%%%%%%% 
%%%%%%%%%%%%%%%%%%%%%%%%%%%%%%%%%%%%%% 
\noindent In this section, we will prove the following theorem:\\

\noindent \textbf{Theorem \ref{thm:main}}\textit{
For all $k \geq 0$, the functors $\sk_k \DDelta$ and $(\sk_0 \DDelta)^{\join (k+1)}$ are both homotopy left cofinal as functors from $\Delta$ to $(k-1)$-connected spaces; in particular, for all %nonempty ?
spaces $X$ and homotopy endofunctors $F$, we have weak equivalences}
\[
\holim_\Delta F(\sk_k \DDelta \join X)  \sim \holim_\Delta  F( (\sk_0 \DDelta)^{\join(k+1)}\join X)
\]
\textit{Furthermore, with Theorem \ref{thm:Tnk}, we have weak equivalences for all $k \geq 0$}
\[
\holim_{\Delta} F(\sk_k \DDelta \join X) %
\sim \holim_n (\cdots \ra \T^{k+1}_n F(X) \ra \T^{k+1}_{n-1} F(X) \ra \cdots \T^{k+1}_1 F(X))
\]
%\bigskip

%
%%%%%%%%%%
% Step one: define intermediate category. 
%%%%%%%%%%
We first make a few definitions. Let $\mathcal{K}$ be the category of $(k-1)$-connected spaces of CW type. Define functors $\mathscr{X}_k, \Y_k : \Delta \ra \mathcal{K}$ %Spaces (actually, we will use use that they land in $(k-1)$-connected spaces of CW type) 
such that
\[
\begin{array}{rcl}
\mathscr{X}_k(p) &=& \sk_k\bigtriangleup^p\\
\mathscr{Y}_k(p)& =& \underset{(k+1) \text{ copies } }{\underbrace{\sk_0\bigtriangleup^p\join \cdots \join \sk_0 \bigtriangleup^p}} 
\end{array}
\]

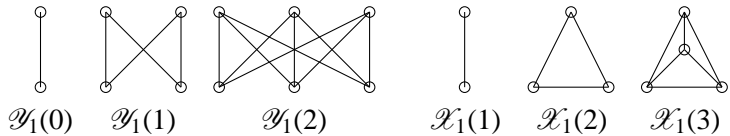
\begin{figure}[h]
\[
\begin{array}{cccccccccc}
%K_1,1
\begin{tikzpicture}
\foreach \x in {0}
\foreach \y in {0}
{\draw (\y,0) -- (\x,1);}
\foreach \x in {0}{
\draw (\x,0) circle (2pt);
\draw (\x,1) circle (2pt);}
\end{tikzpicture}
&
%K_2,2
\begin{tikzpicture}
\foreach \x in {0,1}
\foreach \y in {0,1}
{\draw (\y,0) -- (\x,1);}
\foreach \x in {0,1}{
\draw (\x,0) circle (2pt);
\draw (\x,1) circle (2pt);}
\end{tikzpicture}
&
%K_3,3
\begin{tikzpicture}
\foreach \x in {0,1,2}
\foreach \y in {0,1,2}
{\draw (\y,0) -- (\x,1);}
\foreach \x in {0,1,2}{
\draw (\x,0) circle (2pt);
\draw (\x,1) circle (2pt);}
\end{tikzpicture}
& &
%K_2
\begin{tikzpicture}
\foreach \x in {0}
\foreach \y in {0}
{\draw (\y,0) -- (\x,1);}
\foreach \x in {0}{
\draw (\x,0) circle (2pt);
\draw (\x,1) circle (2pt);}
\end{tikzpicture}
&
%K_3
\begin{tikzpicture}
\foreach \x in {0,1}
\foreach \y in {1}
{\draw (.5,\y) -- (\x,0);}
\foreach \x in {0,1}{
\draw (\x,0) circle (2pt);}
\draw (0,0) --(1,0);
\draw (.5,1) circle (2pt);
\end{tikzpicture}
&
%K_4
\begin{tikzpicture}
\foreach \x in {0,1}
\foreach \y in {.5,1}
{\draw (.5,\y) -- (\x,0);}
\foreach \x in {0,1}{
\draw (\x,0) circle (2pt);}
\draw (0,0) --(1,0);
\draw (.5,.5)--(.5,1);
\foreach \y in {.5,1}{
\draw (.5,\y) circle (2pt);}
\end{tikzpicture}
 \\
\Y_1(0) & \Y_1(1)& \Y_1(2)& & \X_1(1) & \X_1(2)& \X_1(3) \\
\end{array}
\]
\caption{$\Y_1$ of 0,1, and 2 and $\X_1$ of 1,2,and 3}
\label{fig:k=1} 
\end{figure}

\noindent Clearly, when $k=0$, these two functors are the same. Notice also that for $k=1$, we have $\X_1(p) = K_p$ (the complete graph on $(p+1)$ vertices) and $\Y_1(p)=K_{p+1,p+1}$ (the complete bipartite graph on two sets of $(p+1)$ vertices). See figure \ref{fig:k=1}.

%%%%%%%%%%
% Step two: why cofinal is what we want. 
%%%%%%%%%%
Once we show that $\X_k$ and $\Y_k$ are homotopy left cofinal as functors $\Delta \ra \mathcal{K}$, Lemma \ref{thm:hirsch} implies that the maps below are weak equivalences:
\[
\holim_\Delta F(\X_k \join X) \overset{\sim}{\la} \holim_{U \in \mathcal{K}} F(U \join X) \overset{\sim}{\ra}\holim_\Delta F(\Y_k \join X) 
\]
That is, we obtain our desired weak equivalence via a zig-zag of weak equivalences. 
\medskip

%%%%%%%%%%
% Step three: how to show cofinal. 
%%%%%%%%%%
To show left cofinality of the two functors, we show that for all $Z \in (k-1)$-connected spaces of CW type, the simplicial sets $p \mapsto \Top_{cts}(\X_k(p), Z)$ and $p \mapsto \Top_{cts}(\Y_k(p), Z)$ are contractible.

Our proofs will make use of the following lemma, which establishes that for $k=0$, both simplicial sets are contractible: 

%%%%%%
\begin{lem}\label{lem:sk_0}
Let $\sk_0 \DDelta$ be the cosimplicial space sending $[n]$ to the 0 skeleton of $\Delta^n$, the topological $n$-simplex. Then for any nonempty space $Z$, the simplicial set $p \mapsto \Top(\sk_0 \vDelta^p, Z)$ is contractible by a contracting homotopy.
\end{lem}
%%%%%%

\begin{proof} 
For each $p$, $\sk_0 \vDelta^p$ is the discrete space with $(p+1)$ points. This allows us to write $\Top (\sk_0 \vDelta^p, Z)$ as $\prod \limits^{p+1} Z$. That is, our simplicial set is of the form

\[
\xymatrix{
Z \ar[r]& Z \x Z \ar@<-5pt>[l] \ar@<5pt>[l] \ar@<4pt>[r] \ar@<-4pt>[r] &\ar@<-8pt>[l] \ar@<8pt>[l] \ar[l] Z \x Z \x Z\cdots
}
\]

The structure maps are 

\[
\left\{
\begin{array}{rcl}
d_i: Z_{n+1} \ra Z_n &=& \text{projection by deleting the $i$th coordinate}\\
 && \text{e.g. } d_i(z_0,\ldots, z_n) = (z_0, \ldots, z_{i-1}, z_{i+1}, \ldots, z_n) = \bar{z}_{\hat{i}}\\
 \\
s_i: Z_{n} \ra Z_{n+1} &=& \text{inclusion and diagonal applied to the $i$th coordinate} \\
     && \text{e.g. } s_i(z_0,\ldots, z_{n-1}) = (z_0, \ldots, z_{i-1}, z_i, z_i, z_{i+1}, \ldots, z_{n-1})\\
\end{array}
\right.
\]
\smallskip

Recall (from, e.g. Dugger's primer \cite{duggerPrimer}) that for $\mathcal{Y}_\pt$ a simplicial set augmented by the map $d_0: \mathcal{Y}_0 \ra \ast$ (i.e., $\ast =:\mathcal{Y}_{-1}$), a (forward) \textit{contracting homotopy} of $\mathcal{Y}_\pt$ is given by a collection of maps $\mathcal{S}: \mathcal{Y}_n \ra \mathcal{Y}_{n+1}$ for $n \geq -1$ such that for each $y \in \mathcal{Y}_n$, one has 

\[
\left\{
\begin{array}{rcll}
d_i (\mathcal{S} y) &=&
\left\{
	\begin{array}{ll}
	\mathcal{S} (d_i y) & \text{ if } 0 \leq i <n\\
	y & \text{ if } i=n\\
	\end{array}	\right.\\ \\
\mathcal{S}(s_i y) &=&s_i(\mathcal{S} y) & \text{ for } 0 \leq i \leq n\\
\end{array}
\right.
\]
\smallskip

First choose a point $v \in Z$. %
We set $\mathcal{S}(\ast)$(in our -1st dimension) to be $v \in Z$.  For $n$-simplices $z$ for $n>-1$, we define $\mathcal{S}(z) :=(z,v)$. That is, if $z=(z_0, \ldots, z_n) \in \prod \limits^{n+1} Z$, then $\mathcal{S}(z) = (z_0, \ldots, z_n, v)$.  This is our desired contracting homotopy %

\end{proof}

%%%%%%%%%%%%%%%%%%%%%%%%%%%%%%%%%%%%%%%%%
%%%%%%%%%%%%%%%%%%%%%%%%%%%%%%%%%%%%%%%%%
%
\subsection{Contractibility of of $p \mapsto \Top (\X_k (p), Z)$.} 
%
%%%%%%%%%%%%%%%%%%%%%%%%%%%%%%%%%%%%%%%%%
%%%%%%%%%%%%%%%%%%%%%%%%%%%%%%%%%%%%%%%%%

The $k$-skeleton of a (co)simplicial object is adjoint to its $k$-coskeleton. With this, we have, for every $j$, the following isomorphisms of sets:
\begin{spacing}{1.25}
\[
\begin{array}{rcl}
\Top(|\sk_k \Delta^j|, Z) & \cong & \text{sSets}(\sk_k \Delta^j, \text{Sing}(Z))\\
%&\cong& \text{sSets}(\text{rest}_k \vDelta^j, \text{rest}_k \text{Sing}(Z))\\
&\cong& \text{sSets} ( \Delta^j, \text{cosk}_k \text{Sing}(Z))\\
&:=& \cosk_k(\text{Sing}(Z))_j
\end{array}
\]\end{spacing}
That is, $p \mapsto \Top (\X_k (p), Z) \cong \cosk_k (\Sing (Z))$.

For $Y_\pt$ a simplicial set, we have that the map $Y_\pt \ra \cosk_k Y_\pt$ is 1-1 and onto for dimensions $\leq k$, which implies that the homotopy groups of the two objects are the same in dimensions $<k$ (this is discussed in Dwyer-Kan\cite[\S 1.2,part (vi)]{dwyerObstr}).
We also have that the homotopy groups of $\cosk_k Y_\pt$ are trivial in dimensions $\geq k$, when $Y_\pt$ is fibrant. 

Singularization produces fibrant simplicial sets, so we know that $\text{Sing}(Z)$ is fibrant and therefore $\pi_i \cosk_k \text{Sing}(Z) \cong 0$ for all $i\geq k$. We assumed that $Z$ was $(k-1)$ connected, so for $i \leq (k-1), \; \pi_i \cosk_k \text{Sing}(Z) \cong \pi_i \cosk_k Z \cong 0$. We have just shown that all of its homotopy groups are trivial, i.e. it is weakly contractible. Its realization is a space of CW type, so by the Whitehead theorem, weakly contractible implies contractible. \qed

%%%%%%%%%%%%%%%%%%%%%%%%%%%%%%%%%%%%%%%%%
%%%%%%%%%%%%%%%%%%%%%%%%%%%%%%%%%%%%%%%%%
%
\subsection{Contractibility of of $p \mapsto \Top (\Y_k (p), Z)$.} 
%
%%%%%%%%%%%%%%%%%%%%%%%%%%%%%%%%%%%%%%%%%
%%%%%%%%%%%%%%%%%%%%%%%%%%%%%%%%%%%%%%%%%
This proof will be by induction. Note that the base case $k=0$ is Lemma  \ref{lem:sk_0}. We then assume that for all $0<k<K, \; \Top (\Y_k (\pt), Z)$ is contractible. 

For the general case, we will express $\Top (\Y_K(\pt), Z)$ as the homotopy pullback of contractible simplicial sets, and conclude that it is contractible. \\

A common model for the join of two spaces is the following pushout: 
%%%%%%%%%%%%%
\[
\xymatrix{
X \x Y \ar[r]\ar[d] &CX \x Y \ar[d]\\
X\ar[r] &X \join Y \\
}
\]
%%%%%%%%%%%%%

Therefore, for all $j \geq 0$, $\Y_K(j) = \sk_0 \Delta^j \join \Y_{K-1}(j)$ is the following pushout of spaces:
%%%%%%%%%%%%%
\[
\xymatrix{
\Y_{K-1}(j) \x \sk_0 \Delta^j \ar[r]\ar[d] & C\Y_{K-1}(j) \x \sk_0 \Delta^j\ar[d]\\
 \Y_{K-1}(j)  \ar[r] & \Y_K(j):= \Y_{K-1}(j) \join \sk_0 \Delta^j\\
}
\]
%%%%%%%%%%%%%
 
We then apply $\Top (\;, Z)$, the hom-set in topological spaces (of CW type), with $Z$ $(k-1)$-connected. This functor takes pushouts to pullbacks (strict), and we have that the resultant square is, for each $j$, a pullback of sets:

\begin{figure}[h]
%%%%%%%%%%%%%
\[
\xymatrix{
\Top (\Y_{K-1}(j)\join \sk_0 \Delta^j, Z) \ar[r] \ar[d] & \Top (C(\Y_{K-1}(j)) \x \sk_0 \Delta^j, Z)\ar[d] \\
\Top(\Y_{K-1}(j), Z) \ar[r] &\Top(\Y_{K-1}(j) \x \sk_0 \Delta^j, Z)
} 
\]
\caption{Levelwise pullback of sets}
\label{fig:Ykpb}
%%%%%%%%%%%%%
\end{figure}

We will first show contractibility of the simplicial sets, then that the square is a homotopy pullback of simplicial sets (when $j$ is allowed to vary).

By the induction hypothesis, we have that $\Top(\Y_{K-1}(\pt), Z)$ is a contractible simplicial set. For the other two simplicial sets,  if we allow our indices to vary independently, we may view our square as one of bisimplicial sets, and make use of results of Bousfield and Friedlander. Note that we consider $\Top(\Y_{K-1}(\pt), Z)$ bisimplicial by making it constant in one direction. See Figure \ref{fig:bi-sq}.

\begin{figure}[h]
%%%%%%%%%%%%%
\[
\xymatrix{
\Top (\Y_{K-1}(i)\join \sk_0 \Delta^j, Z) \ar[r] \ar[d] & \Top (C(\Y_{K-1}(i) \x \sk_0 \Delta^j, Z)\ar[d] \\
\Top(\Y_{K-1}(i), Z) \ar[r] &\Top(\Y_{K-1}(i) \x \sk_0 \Delta^j, Z)
} 
\]
\caption{Square of bisimplicial sets (the indices $i,j \geq 0$ vary independently)}
\label{fig:bi-sq}
%%%%%%%%%%%%%
\end{figure}

Using the following result, we can show levelwise equivalence (i.e. contractibility) and conclude that the diagonal is also contractible:

\begin{thm}(Bousfield-Friedlander \cite[p119, Theorem B.2]{bousfield-friedlander})\label{thm:bf2}	
Let $f: X\ra Y$ be a map of bisimplicial sets such that $f_{m,\ast}: X_{m,\ast}\ra Y_{m,\ast}$  is a weak equivalence for each $m \geq 0$. Then $diag(f): diag X\ra diag Y$ is a weak equivalence.
\end{thm}

\noindent Note that we have for each $m\geq 0$ the following isomorphisms: 
\[
\begin{array}{lcl}
\Top (C(\Y_{K-1}(m)) \x \sk_0 \DDelta, Z) &\cong& \Top (\sk_0 \DDelta, \hom_{\Top}(C(\Y_{K-1}(m)), Z)) \\
\Top (\Y_{K-1}(m) \x \sk_0 \DDelta, Z) &\cong &\Top (\sk_0 \DDelta, \hom_{\Top}(\Y_{K-1}(m) , Z))\\
\end{array}
\]
Then, by Lemma \ref{lem:sk_0}, these are contractible, since we have expressed them as $\Top( \sk_0 \DDelta, X)$ for $X$ a space. By Theorem \ref{thm:bf2}, we may also conclude that $\Top (C(\Y_{K-1}(\pt)) \x \sk_0 \DDelta, Z)$ and  $\Top (\Y_{K-1}(\pt) \x \sk_0 \DDelta, Z)$ are contractible, as we have shown levelwise contractibility. %, therefore they are connected. 

To show that Figure \ref{fig:Ykpb} is a \textit{homotopy} pullback, we continue to regard our square as one of bisimplicial sets, as in Figure \ref{fig:bi-sq}, and use the following Theorem:

\begin{thm}(Bousfield-Friedlander \cite[Theorem B.4]{bousfield-friedlander})\label{thm:bf}
Let 
\[
\xymatrix{
V \ar[r]\ar[d] & X\ar[d]\\
W \ar[r] & Y\\
}
\]
be a commutative square of bisimplicial sets such that the terms $V_{m, \ast}, W_{m, \ast}, X_{m, \ast}$ and $Y_{m, \ast}$ form a homotopy fiber square for each $m\geq 0$. If $X_{m, \ast}$ and $Y_{m, \ast}$ are all connected, then 
\[
\xymatrix{
\mathrm{diag}V \ar[r]\ar[d] & \mathrm{diag} X\ar[d]\\
\mathrm{diag}W\ar[r]& \mathrm{diag}Y\\
}
\]
is a homotopy fiber square.
\end{thm}

First note that $X$ and $Y$ in our case have already been shown to be contractible, so are therefore connected.  We will show that our bisimplicial set diagram levelwise is homotopy pullback squares of simplicial sets. That is, for all $m \geq 0$, the following is not just a levelwise pullback of sets but a homotopy pullback of simplicial sets: 

\[
\xymatrix{
\Top (\Y_{K-1}(m)\join \sk_0 \DDelta, Z) \ar[r] \ar[d] & \Top (C(\Y_{K-1}(m)) \x \sk_0 \DDelta, Z)\ar[d] \\
\Top(\Y_{K-1}(m), Z) \ar[r] &\Top(\Y_{K-1}(m) \x \sk_0 \DDelta, Z)
} 
\]
We then show that the righthand vertical map is a Kan fibration\footnote{We only need to show one map is a fibration due to right properness of the category of simplicial sets. See the glueing lemma, e.g. in Schwede\cite{schwede97}, Lemma 1.19.}, and conclude that our square is a homotopy pullback.

Simplicial sets are a simplicial model category, satisfying Quillen's SM7 axiom, which is as follows\footnote{$\Hom_C( , )$ is the simplicial set from the simplicial model structure,  $\Hom_C(X,Y) := C(X \x \DDelta, Y)$ where $C(X,Y)$ is the hom-set of morphisms in $C$.}:  \\

\noindent \textbf{Quillen's SM7 axiom}: \textit{
Let $\C$ be a simplicial model category, and $\Hom_\C(X,Z)$ denote the simplicial set of morphisms between $X$ and $Z$. If $Y\in \C$ is fibrant and $f: A\ra B$ is a cofibration in $\C$, then $\Hom_\C (f, Y): \Hom_\C (B, Y) \ra \Hom_\C (A, Y)$ 
is a fibration of simplicial sets. 
}\\

We apply several adjunctions to get the following isomorphisms. Note that $\sk_0 \DDelta$ is equivalent to the levelwise simplicial set skeleton, i.e. $\sk_0 \DDelta: n \mapsto \sk_0 \Delta^n$, so it is adjoint levelwise to the simplicial coskeleton (for discussions of skeleta/coskeleta, see Goerss-Jardine \cite[Ch. VII, mainly \S 1]{GJ}). %We will denote the singularization of a space $Z$ by $sZ$.

\begin{spacing}{1.25}
\[
\begin{array}{lcl}
 \Top (C \Y_{K-1}(m) \x \sk_0 \DDelta, Z) &\cong& \sSet (C \Y_{K-1}(m) \x \sk_0 \DDelta,\text{Sing}(Z) ) \\
 								&\cong&\sSet (C\Y_{K-1}(m), \Hom_{sSet} (\sk_0 \DDelta, Z))\\
 								&\cong&\sSet (C\Y_{K-1}(m), \Hom_{sSet} (\DDelta, \cosk_0 \text{Sing}(Z)))\\
								&\cong& \sSet (C\Y_{K-1}(m) \x \DDelta, \cosk_0 \text{Sing}(Z))\\
								&=:& \Hom_{sSet}(C\Y_{K-1}(m), \cosk_0 \text{Sing}(Z))\\
															%	\\
								\text{and likewise, }\\
								\\
\Top(\Y_{K-1}(m) \x \sk_0 \DDelta, Z)	&\cong&\sSet (\Y_{K-1}(m) \x \sk_0 \DDelta, \text{Sing}(Z))	\\
 								&\cong& \sSet (	\Y_{K-1}(m) \x \DDelta, \cosk_0 \text{Sing}(Z)) \\
								&=:& \Hom_{sSet}(\Y_{K-1}(m), \cosk_0 \text{Sing}(Z)) \\ 						
 \end{array}
\]
\end{spacing} 

We may express our righthand vertical map as 
\[
\Hom_{sSet} (C\Y_{K-1}(m), \cosk_0 \text{Sing}(Z)) \ra  \Hom_{sSet}(\Y_{K-1}(m), \cosk_0 \text{Sing}(Z)) 
\]
which is $\Hom_{sSet}( \; \;, \cosk_0 \text{Sing}(Z))$ applied to the map $\Y_{K-1}(m) \ra C\Y_{K-1}(m)$. This map is a cofibration of simplicial sets since it is a monomorphism. The singularization of a topological space is a fibrant simplicial set, and as coskeleton is a right Kan extension, it preserves fibrant objects. Therefore, we may apply SM7 and conclude that our map is a fibration of simplicial sets, and our square is a homotopy pullback square. Applying Theorem \ref{thm:bf}, we conclude that Figure \ref{fig:Ykpb} is a homotopy pullback of simplicial sets which we have already shown to be contractible, so we conclude that $\Top (\Y_{K-1}(\pt)\join \sk_0 \Delta^\pt  Z)$ is also a contractible simplicial set. \qed

%\newpage
%%%%%%%%%%%%%%%%%%%%%%%%%%%%%%%%%%%%
%%%%%%%%%%%%%%%%%%%%%%%%%%%%%%%%%%%%
%%%%%%%%%%%%%%%%%%%%%%%%%%%%%%%%%%%%
\bibliographystyle{gtart}
%\bibliographystyle{amsalpha}
%\bibliography{thesisbib}
%%%%%%%%%%%%%%%%%%%%%%%%%%%%%%%%%%%%
%%%%%%%%%%%%%%%%%%%%%%%%%%%%%%%%%%%%
%%%%%%%%%%%%%%%%%%%%%%%%%%%%%%%%%%%%
\providecommand{\bysame}{\leavevmode\hbox to3em{\hrulefill}\thinspace}
\providecommand{\MR}{\relax\ifhmode\unskip\space\fi MR }
\providecommand{\MRhref}[2]{%
  \href{http://www.ams.org/mathscinet-getitem?mr=#1}{#2}
}
\providecommand{\href}[2]{#2}

\end{document}